\documentclass[11pt,a4paper]{article}
\usepackage{ifthen}
\usepackage{amsmath}
\usepackage{amssymb}
\usepackage{latexsym}
\usepackage{bbm}
\usepackage{graphicx}
\usepackage{pdfsync}
\usepackage{verbatim}
\usepackage{fancyhdr}
\usepackage{url}               
\usepackage[english]{babel}
\usepackage[applemac]{inputenc}

\usepackage[pdftex,%
a4paper=true,%
backref=true,%
pagebackref=true,%
bookmarks=true,%
bookmarksnumbered=true,%
colorlinks=true,%
urlcolor=blue,%
pdfstartview=FitH%
]{hyperref}

\newtheorem{theorem}{Theorem}[section]
\newtheorem{lemma}[theorem]{Lemma}

\newtheorem{proposition}[theorem]{Proposition}
\newtheorem{corollary}[theorem]{Corollary}

\newtheorem{remark}[theorem]{\it \bf{Remark}\/}

\numberwithin{equation}{section}
\numberwithin{figure}{section}

\newcommand{\bp}{{\it Proof. }}
\newcommand{\ep}{\hfill $\square$\\}

\newcommand{\be}{\begin{equation}}
\newcommand{\ee}{\end{equation}}
\newcommand{\bea}{\begin{eqnarray}}
\newcommand{\eea}{\end{eqnarray}}
\newcommand{\bee}{\begin{eqnarray*}}
\newcommand{\eee}{\end{eqnarray*}}

\def\CC{\mathbb{C}}
\def\NN{\mathbb{N}}
\def\RR{\mathbb{R}}

\def\dps{\displaystyle}

\def\eps{\vare}

\def\eps{\varepsilon}
\def\LL{\mathcal{L}}

\catcode`@=11
\def\supess{\mathop{\operator@font Sup\,ess}}
\catcode`@=12
\def\CC{\mathbb{C}}

\def\NN{\mathbb{N}}

\def\RR{\mathbb{R}}
\def\CC{\mathbb{C}}
\def\LL{\mathbb{L}}

\def\KK{\mathbb{K}}

\def\R2+{\RR ^2_+}

\def\lim{\mathop{\rm lim}}

\begin{document}

\renewcommand{\refname}{References}
\bibliographystyle{alpha}

\pagestyle{fancy}
\fancyhead[L]{ }
\fancyhead[R]{}
\fancyfoot[C]{}
\fancyfoot[L]{ }
\fancyfoot[R]{}
\renewcommand{\headrulewidth}{0pt} 
\renewcommand{\footrulewidth}{0pt}

\newcommand{\montitre}{Diffusion approximation for Fokker Planck with heavy tail equilibria : a spectral method in dimension 1 }

\newcommand{\auteur}{\textsc{ Gilles Lebeau, Marjolaine Puel}}
\newcommand{\affiliation}{Laboratoire J.-A. Dieudonn\'e \\
 Universit\'e de Nice Sophia-Antipolis\\
 Parc Valrose, 06108 Nice Cedex 02,   France\\
\url{lebeau@unice.fr, Marjolaine.Puel@unice.fr}
}

 \begin{center}
{\bf  {\LARGE \montitre}}\\ \bigskip \bigskip
 {\large\auteur}\\ \bigskip \smallskip
 \affiliation \\ \bigskip
\today
 \end{center}

\begin{abstract}
This paper is devoted to the diffusion approximation for the $1$-d Fokker Planck equation with 
a heavy tail equilibria of the form $(1+v^2)^{-\beta/2}$, in the range $\beta\in ]1,5[$. We prove that the 
limit  diffusion equation involves a fractional Laplacian $\kappa\vert \Delta\vert^{\frac{\beta+1}{6}}$,
and we compute the value of the diffusion coefficient $\kappa$. This extends previous results of 
E. Nasreddine and M. Puel  \cite{NaPu} in the case $\beta>5$, and of 
P. Cattiaux, E. Nasreddine and M. Puel  \cite{CaNaPu} in the case $\beta=5$. 

\footnote{ Gilles Lebeau  was supported by the European Research Council, 
ERC-2012-ADG, project number 320845:  Semi Classical Analysis of Partial Differential
Equations.}
\end{abstract}

\tableofcontents
\newpage

 \pagestyle{fancy}
\fancyhead[R]{\thepage}
\fancyfoot[C]{}
\fancyfoot[L]{}
\fancyfoot[R]{}
\renewcommand{\headrulewidth}{0.2pt} 
\renewcommand{\footrulewidth}{0pt} 

\section{Introduction}

\subsection{Setting of the problem}

  In this present paper, we deal with the equation
\begin{equation}
 \partial_t f+ v\cdot \nabla_x f=Q(f)
\end{equation}
where the the Fokker Plank operator $Q$ is given by 
\begin{equation}\label{defQ}
Q(f)=\nabla_v\cdot (\frac{1}{\omega}\nabla_v (\omega f))
\end{equation}
for a fixed $\omega$ that determines the equilibrium $F=\frac{C_\beta^2}{\omega}$, $C_\beta$ beeing a normalization constant.

Recall that the aim of diffusion approximation is to provide a simpler model when the interaction between particles are the dominant phenomena and when the observation time is very large. For that purpose, we introduce a small parameter, $\eps$, the mean free path and we proceed to a rescaling in time and space
$$
t=\frac{t'}{\theta(\eps)}\quad x=\frac{x'}{\eps}
$$
which leads to the following rescaled equation (without primes)
\begin{equation}\label{fp-theta}
\theta(\eps) \partial_t f^\eps+\eps v\cdot \nabla_x f^\eps=Q(f^\eps).
\end{equation}

\noindent Passing formally to the limit, we get that 
$$
f^\eps\rightarrow_{\eps\rightarrow0}f^0=\rho(t,x) F(v)
$$
where $F(v)$ is the equilibrium defined above. It remains to identify the equation satisfied by $\rho$.

When the equilibrium $F$ is a gaussian, it is classical (see \cite{BaSaSe},\cite{BeLi},\cite{DeGoPo},\cite{LaKe},\cite{De2} for Boltzmann and \cite{DeMa-Ga} for Fokker Planck) that by taking the classical time scaling $\theta(\eps)=\eps^2$, we obtain for $\rho$ a diffusion equation
\begin{equation}\label{d}
\partial_t\rho-\nabla_x(D\nabla\rho)=0
\end{equation}
where 
\begin{equation}\label{dc}
D=\int vQ^{-1}(-vF) dv.
\end{equation}
Indeed, the formal expansion $f^\eps=f^0+\eps f^1+\eps^2 f^2...$gives
$$
\begin{array}{rcl}
Q(f^0)&=&0\\
Q(f^1)&=& v\cdot \nabla_x f^0\\
Q(f^{2,\eps})&=&\partial_tf^0+v \cdot \nabla f^1
\end{array}
$$
and the compatibility equation for the equation giving $f^2$ gives 
$$
\partial_t \int f^0+ \int \nabla \cdot (v Q^{-1} (v\cdot \nabla_x f^0))=0
$$
which is another formulation of (\ref{d}) since $f^0=\rho(t,x) F(v)$ and $F$ is normalized by $\dps \int F=1$.

In the present work, we consider heavy tail equilibria $F(v)=\frac{C_\beta^2}{\omega}$ with $\omega=(1+|v|^2)^\frac{\beta}{2}$. In \cite{NaPu}, the classical scaling is studied and it is proved that we obtain a diffusion equation (\ref{d}), (\ref{dc}) as soon as $\beta>d+4$. The critical case where $\beta=d+4$ is studied in \cite{CaNaPu} where the expected result of classical diffusion with an anomalous time scaling is proved.

The aim of this paper is to study the case where $\beta<d+4$, when the diffusion  coefficient (\ref{dc}) is not defined anymore.
We need to operate an ad hoc rescaling in time  that we will compute during the proof. Fractional diffusion limit has been already obtained in the case of the linear Boltzmann equation for heavy tail equilibria when the cross section is such that the operator has a spectral gap (see \cite{MMM} for the pioneer paper in the case of space independent cross section, \cite{Me} for a weak convergence result and \cite{NBAMePu1} for a strong convergence result) and when the cross section is degenerated \cite{NBAMePu1}. The main difficulty of this case is due to the fact that the Fokker Planck operator $Q$ has no spectral gap. The idea here is thus to study the whole operator, advection plus collision, at $\eps$ fixed to compute the first eigenvalue and its corresponding eigenvector. The dependency of the first eigenvalue with respect to $\eps$ will give us the right time scaling and the power of the limiting fractional diffusion operator. Note that a fractional diffusion has also been obtained for a Fokker-Planck like operator in \cite{CeMeTr}.
\medskip

\noindent{\bf Outline of the paper}

In the next subsection, we recall the previous results obtained for this equation with heavy tail equilibria, and we quote the main theorem of this present paper and proceed to a change of unknown. It is followed by a section dedicated to the computation of the first eigenfunction and eigenvalue. Finally, in section \ref{secmoment}, we apply the momentum method to complete the proof of the main theorem.

\subsection{Previous results}
 The functional setting of the study of equation \ref{fp-theta} has been settled in \cite{NaPu} where we define the functional ad hoc spaces $Y^p_\omega\left(\mathbb{R}^{2d}\right)=L^p\left(\mathbb{R}^d,\  H_p(\mathbb{R}^d)\right)$, where
\begin{equation}H_p(\mathbb{R}^d)= \left\{f: \mathbb{R}^d\rightarrow \mathbb{R},\int_{\mathbb{R}^{d}} |f|^p \ \omega^{p-1}\ \mathrm{d}v< \infty\right\},\end{equation}  where $\omega=(1+||v||^2)^{\frac{\beta}{2}}$ and
$$
L^\infty_\omega(\RR^d)=\{f: \mathbb{R}^d\rightarrow \mathbb{R}, f\omega \in L^\infty(\RR^d)\}.
$$

Define
\begin{equation}\label{V}
V=\left\{f: \mathbb{R}^d\rightarrow \mathbb{R}, \int_{\mathbb{R}^d} |f|^2 \ \omega\ \mathrm{d}v< \infty\ \mathrm{and}\ \int_{\mathbb{R}^d} \frac{|\nabla_v(f\ \omega)|^2}{\omega}\ \mathrm{d}v <\infty\right\},
\end{equation}
 $V'$ being its dual.
 
 \medskip
 
\noindent{\bf{Operator's properties.}}  We sumerize in the following proposition the main properties of the interaction operator.
\begin{proposition}\label{qf}\cite{NaPu}
Let $f$ and $g$ be smooth functions in $V$ defined in \eqref{V}. The following assertions hold true:
\begin{enumerate}
\item The operator $Q$ is conservative, thus equation (\ref{fp-theta}) preserves the total mass of the distribution
$$\int_{\mathbb{R}^d} Q(f)\ dv=0,\ \ \mbox{for all  } f\in V.$$
\item The operator $Q$ is self-adjoint with respect to the measure $\omega\ \mathrm{d}v$:
\begin{equation} \int_{\mathbb{R}^d} Q(f)\ g\  \omega\ dv=-\int_{\mathbb{R}^d} \frac{\nabla_v(f\ \omega)\cdot \nabla_v(g\ \omega)}{\omega}\ dv=\int_{\mathbb{R}^d} f\ Q(g)\  \omega\ dv,\end{equation}
\item The operator $Q$ is dissipative:
\begin{equation}\label{prod} \int_{\mathbb{R}^d} Q(f)\ f\  \omega\ dv =-\int_{\mathbb{R}^d} \frac{|\nabla_v(f\ \omega)|^2}{\omega}\ dv\leq 0.\end{equation}
 \item The kernel of $Q$ is one-dimensional and spanned by $\frac{1}{\omega}$.
 \item   The operator $Q$ is continuous from $V\longrightarrow V'.$
 \end{enumerate}\end{proposition}
 
 \medskip
 
\noindent{\bf{Existence theorem.}} We recall  the following theorem inspired from \cite{De1}
\begin{theorem}\label{ex} \cite{NaPu}Let $\eps$ be fixed. Assume that $f_0\in Y^2_\omega(\mathbb{R}^d)$, equation \eqref{fp-theta} has a unique solution $f$ in the class of functions $Y$ defined by:
$$Y= \left\{f \in L^2\left([0, T]\times \mathbb{R}^d, V\right),\ \ \theta(\eps)\partial_t f+\eps v\cdot\nabla _x f\in L^2\left([0, T]\times \mathbb{R}^d,V'\right)\right\}.$$
\end{theorem}

 \medskip
 
\noindent{\bf{Classical diffusion approximation.}} The case where $\beta>d+4$ leads to a diffusion equation as described in the following theorem.
\begin{theorem}\label{th}\cite{NaPu}
Assume now that $\beta>d+4$. Assume that  $f_0$ is a nonnegative function in $ Y^2_\omega\cap Y^p_\omega$ with $p>2$. Assume that $\theta(\eps)=\eps^2$, let $f^\varepsilon$ be the solution of \eqref{fp-theta} in $Y$ with initial data $f_0$,.\\
Then, $f^\varepsilon$ converges weakly star in $L^\infty\left([0, T],\  Y^p_\omega(\mathbb{R}^{2d})\right)$ towards $\rho(t,x)\ \frac{C_\beta^2}{\omega}$ where $\rho(t,x)$ is the unique solution of the system
\begin{equation}\label{diff}
\partial_t\rho+\nabla_x\cdot j=0
 \end{equation}
 \begin{equation}\label{J}
 j=-D\ \nabla_x \rho,
\end{equation}
where the initial datum is given by $\dps\rho_0(x)=\int_{\mathbb{R}^d}\ f_0\ dv,$ and the diffusion tensor $D$ is given by
\begin{equation}\label{dii}
D=\int_{\mathbb{R}^d}\ v\otimes\chi\ dv,
\end{equation}
where $\chi$ is the unique solution of the cell equation $Q(\chi)=\frac{-C_\beta^2\ v}{\omega}$ with $\dps\int_{\mathbb{R}^d}\ \chi\ dv=0$.
\end{theorem}

 \medskip
 
\noindent{\bf{Critical case, $\beta=d+4$}.}
\begin{theorem}\label{thmcritik}\cite{CaNaPu}
Assume that $\beta=d+4$. Then there exists $\kappa>0$ such that, with $\theta(\varepsilon)=\varepsilon^2 \, \ln(1/\varepsilon)$, for all initial density of probability $f_0$, the solution $f^\varepsilon_t$ of \eqref{fp-theta} weakly converges as $\varepsilon \to 0$ towards $$(v,x) \mapsto C_\beta^2 \, \omega_\beta^{-1}(v) \, (h_0*\rho_t)(x)$$ where $\rho_t$ is the density of a centered gaussian random vector with covariance matrix $(2\kappa/3) \, t \, \text{Id}$ and $h_0(x)=\dps \int f_0(x,v) \, dv$.
\end{theorem}

\subsection{Main theorem}
Assume from now on that the dimension $d=1$.
\label{secmain}
\begin{theorem}\label{main}
Assume now that $1<\beta<5$ with $\beta\neq \{2,3,4\}$. Assume that  $f_0\in L^1$ is a nonnegative function in 
$ Y^2_\omega$ and $f_0 \omega\in L^\infty$. Let $f^\varepsilon$ be the solution of \eqref{fp-theta} in $Y$ with initial data $f_0$, when $\theta(\varepsilon)=\varepsilon^{\frac{\beta+1}{3}} .$\\
Let $\kappa=2C_\beta^2 (\beta+1)9^{-\frac{\beta+1}{3}} 
\cos(\frac{\pi}{2}\frac{\beta+1}{3})\frac{\Gamma(1-\frac{\beta+1}{3})}{\Gamma(1+\frac{\beta+1}{3})} >0$, where $\Gamma$ is the Euler function.\\
Then $f^\varepsilon$ converges weakly star in $L^\infty\left([0, T],\  Y^2_\omega(\mathbb{R}^{2})\right)$ towards $\rho(t,x)\ \frac{C_\beta^2}{\omega}$ where $\rho(t,x)$ is the inverse Fourier transform of the unique solution 
$ \hat\rho(t,k)=\int e^{-ixk} \hat\rho(t,x)dx$ of
\begin{equation}\label{diff}
\partial_t\hat \rho+\kappa|k|^{\frac{\beta+1}{3}}\hat\rho=0, \quad \rho(0)=\int f_0 dv,  \ .
 \end{equation}
\end{theorem}

\begin{remark}The hypothesis $\beta\neq \{2,3,4\}$ is  technical. It avoids to introduce logarithmic terms
in the expression of the solution $g$ in \eqref{eqdeg} .  Observe that $\mu=\frac{\beta+1}{3}\in ]2/3,2[$, and that for $\mu\in ]2/3,2[$, one has
$f(\mu)=\cos(\frac{\pi}{2}\mu)\frac{\Gamma(1-\mu)}{\Gamma(1+\mu)}>0$, $f(1)=\pi/2$ and
$\lim_{\mu\rightarrow 2} f(\mu)=+\infty$. \end{remark}

As we said, in order to prove this theorem, we compute the first eigenvalue and eigenvector of the whole operator $(-i\eps v\cdot \nabla+Q)$ and for that purpose, to simplify the computation, we proceed to a change of unknown such that the new operator splits into a Schr\"odinger operator.

 \medskip
 
\noindent{\bf{Changing the unknown.}} We start with the Fokker Planck equation
$$
\partial_tf+v\cdot\nabla_xf=Q(f)=\nabla_v(F\nabla_v(\frac{f}{F}))
$$
with equilibria given by
$$
F=\frac{C_\beta^2}{(1+|v|^2)^{\frac{\beta}{2}}}=\frac{C_\beta^2}{(1+|v|^2)^{\gamma}}
$$
Since we impose $\gamma=\frac{\beta}{2}>\frac{1}{2}$, $F\in L^1(\RR)$,  and we chose $C_\beta$ such that $\dps\int F dv=1$.
In order to work with a self adjoint operator in $L^2$, we proceed to a change of unknown by writing 
$$
f=F^\frac{1}{2} g
$$
and the equation becomes
$$
\partial_t g+v\cdot \nabla_x g=F^{-\frac{1}{2}}\nabla_v(F\nabla_v(\frac{g}{F^{\frac{1}{2}}}))
$$
that can be written 
$$
\partial_t g+v\cdot\nabla_x g=\Delta_v g-W(v) g
$$
with 
$$W(v)=-\frac{1}{2}F^{-\frac{1}{2}}\nabla\cdot(F^{-\frac{1}{2}}\nabla F).
$$
The explicit formula for $W$ is
$$
W(v)= \frac{\gamma}{(1+|v|^2)^2}[|v|^2(\gamma+1)-1]
$$
and its asymptotic behavior for high velocities is 
$$
W(v)\sim_{v\sim\infty}\frac{\gamma(\gamma+1)}{|v|^2}.
$$
We see the equation as 
$$
\partial_tg=-\mathcal{L}g
$$
where $\mathcal{L} =-\Delta_v+W(v)+v\cdot\nabla_x$ is a non negative  operator since
$$
(\mathcal{L}g|g)=\int|\nabla g|^2+\int W(v)|g|^2=\int F\vert \nabla_v(\frac{g}{F^{1/2}}) \vert^2\geq 0\ ,
$$
thus
$$
g=e^{-t\mathcal{L}}g_0.
$$
Since the operator has coefficient that do not depend on $x$, we operate a Fourier transform in $x$ and proceed to a  second change of unknown by writing
$$
g(s,x,v)=(2\pi)^{-1}\int e^{ix\cdot \xi} \tilde g(s,\xi,v) d\xi
$$
where $\tilde g $ satisfies
$$
\partial_t \tilde g=-\mathcal{L}\tilde g
$$
where 
$$
\mathcal{L}\tilde g=-\Delta_v\tilde g+W(v)\tilde g+i(\xi\cdot v) \tilde g
$$

 \medskip
 
\noindent{\bf{Rescaling.}} We do a rescaling both in space and time
$$
t=Ts,\quad x=T^{1-\delta}y,\quad \xi=T^{\delta-1}k
$$
so that $e^{ix\cdot\xi}=e^{iy\cdot k}$. The equation becomes
\begin{equation}\label{rescaled}
\partial_s\tilde g=-T\mathcal{L}_\eps (\tilde g)
\end{equation}
with $\eps=T^{\delta-1}$ and
$$
\mathcal{L}_\eps (\tilde g^\eps)=-\Delta_v \tilde g^\eps+W(v)\tilde g^\eps+i(v\cdot \eps k)\tilde g^\eps \ .
$$

\noindent Classical diffusion corresponds to $\delta=\frac{1}{2}$. When $\beta<5$, the right scaling will be given by the power of $\eps$ of the leading term of the first eigenvalue of the full operator.

\section{Spectral study of the operator : }

In this section, for $\eps >0$, we compute the eigenvalue $\mu^\eps$ with lowest absolute value  and the associated eigenfunction
$M^\eps$ (normalized by $M^\eps(0)$=1) of the unbounded operator $\mathcal{L}_\eps$ acting on $L^2$:
\begin{equation}\label{eqmeps}
\mathcal{L}_\eps M^\eps= -\Delta_v M^\eps+W(v) M^\eps+i (v\cdot \eps k)M^\eps =\mu^\eps M^\eps ,\quad M^\eps(0)=1,\mbox{ and } {M^\eps}'(0)=b.
\end{equation}
In dimension 1, the equation leading to the eigenvalue can be written 
$$
[-\partial^2_v+W(v)+i\eps k v-\mu^\eps]M^\eps=0
$$
and $W$ is given by
$$
W(v)= \frac{\gamma}{(1+|v|^2)^2}[v^2(\gamma+1)-1].
$$
The domain of $\mathcal{L}_\eps$ is 
$$
D(\mathcal{L}_\eps)=\{ g\in L^2,\quad \partial^2_v g \in L^2,\quad vg\in L^2\} \ .
$$
Note that for $\eps >0$, the domain of $\mathcal{L}_\eps$ is not equal to the domain of the limiting operator. 

\noindent In dimension 1, the domain is compact, thus the spectrum is discrete. The construction of the eigenvalue turns out to be a connexion problem between
$$
E^{\pm}_\mu=\{g|\mathcal{L}_\eps(g)-\mu g=0 \quad \mbox{with } g\in L^2(\pm v\geq 0)\}.
$$

\subsection{Large velocities asymptotic : solution to an approximated equation}

Since $W(v)\sim_{|v|\rightarrow \infty} \frac{\gamma(\gamma+1)}{v^2}$, we will first consider the approximated differential equation 
$$
(-\partial^2_v+\frac{\gamma(\gamma+1)}{v^2}+i\eps k v)f=\mu f, \quad v\in ]0,\infty[\, 
$$
or
$$
-v^2\partial^2_vf+\gamma(\gamma+1)f+i\eps k v^3f=\mu v^2 f.
$$
If we want to get rid of the parameter $\eps$, we need to proceed to the following rescaling
$$
v=(\eps k)^{-\frac{1}{3}}s,\quad \mu= (\eps k)^{\frac{2}{3}}\lambda
$$
that leads to 
\begin{equation}\label{eqdiffsing}
(-s^2\partial^2_s +\gamma(\gamma+1)+is^3)f=\lambda s^2 f, \quad s\in ]0,\infty[\ .
\end{equation}
Near $s=0$, equation (\ref{eqdiffsing}) is a differential equation with regular singular points. We proceed to a change of unkown by writing $f=s^\delta g$, with $\delta(\delta-1)=\gamma(\gamma+1)$,  i.e $\delta =-\gamma$ or $\delta = \gamma+1$. Then the new unknown $g$ satisfies 
\begin{equation}\label{eqdeg}
-\partial_s^2g-\frac{2\delta}{s}\partial_s g +(is-\lambda)g=0.
\end{equation}
Writing $\dps g=\sum_0^\infty g_n s^n$  leads to the following equation for $g_n$
$$
-\sum_{n\geq 0}(n+2)(n+1) g_{n+2}s^n-2\delta \sum_{n\geq 0}(n+2) g_{n+2}s^n-2\delta \frac{g_1}{s}-\lambda \sum_{n\geq 0 } g_n s^n+i\sum_{n=1}^\infty g_{n-1}s^n=0\ ,
$$
that gives assuming $\gamma\neq \frac{k+1}{2}, k\in \mathbb N$,
\begin{equation}\label{lesgn}
\left\{\begin{array}{rcl}
g_1&=&0\\
g_{n+2}&=&\frac{1}{(n+2)(n+1+2\delta)}[-\lambda g_n+ig_{n-1}]\quad \forall n\geq 0\quad (g_{-1}=0)
\end{array}\right.
\end{equation}
that define a unique solution if $g_0=1$.
Define 
\begin{equation}\label{defdesF}\begin{array}{l}
F_{+,\lambda}(s)= \sum_0^\infty g_n s^n,\quad g_0=1, g_1=0, g_n \mbox{ defined as above with } \delta=-\gamma\\
F_{-,\lambda}(s)= \sum_0^\infty g_n s^n,\quad g_0=1, g_1=0, g_n \mbox{ defined as above with } \delta=\gamma+1.
\end{array}
\end{equation}
A basis of the solution space of equation (\ref{eqdiffsing}) is thus given by the two independent solutions
$$
F_{+,\lambda}(s) s^{-\gamma}\quad \mbox{and}\quad F_{-,\lambda}(s) s^{\gamma+1}\ .
$$
$F_{\pm,\lambda}$ are normalized by $F_{\pm,\lambda}(0)=1$ and are entire functions of $s\in\mathbb C.
$

\begin{proposition}\label{propH}
Let $F_{+,\lambda}$ and $F_{-,\lambda}$ be defined in (\ref{defdesF}). There exists $\lambda_0$ such that 
for all $ \lambda\in \CC$, such that $|\lambda|\leq \lambda_0$, equation (\ref{eqdiffsing}) has a unique solution $H_\lambda(s)$ such that
\begin{enumerate}
\item $\dps\int_1^\infty |H_\lambda(s)|^2ds<\infty$.

\item $\dps H_\lambda(s)= s^{-\gamma}F_{+,\lambda}(s)+d(\lambda)F_{-,\lambda}(s)s^{\gamma+1}$.
\end{enumerate}
\end{proposition}

\bp
For $s>>1$, we first consider the approximate equation 
\begin{equation}
(-\partial_s^2+is-\lambda)g=0
\end{equation}
that define a unique $L^2(s>R)$ solution (up to a constant) given by 
$$G_\lambda(s)= Ai[e^{i\frac{\pi}{6}}(s+i\lambda)]$$
where $Ai$ is the Airy function given by
$$
Ai(z)=\frac{1}{2\pi}\int_{-\infty}^{\infty}e^{i(\frac{t^3}{3}+zt)}dt.
$$
We will look at a solution to (\ref{eqdiffsing}) via the following change of unknown 
$$H_\lambda =CAi[e^{i\frac{\pi}{6}}(s+i\lambda)][1+R_\lambda(s)] $$
where $R_\lambda$ satisfies
$$
2G'_\lambda R'_\lambda +G_\lambda R''_\lambda=\frac{\gamma(\gamma+1)}{s^2}G_\lambda(1+R_\lambda)
$$
that can be written
$$
R''_\lambda+2\frac{G'\lambda}{G_\lambda}R'_\lambda=\frac{\gamma(\gamma+1)}{s^2}(1+R_\lambda)
$$
and that leads to the implicit equation
$$
R_\lambda(s)=\int_s^\infty [\int_s^z\frac{G_\lambda^2(z)}{G^2_\lambda(u)}du]\frac{\gamma(\gamma+1)}{z^2}(1+R_\lambda(z))dz.
$$
We need to prove the following lemma.
\begin{lemma}\label{deuxdeux} Define 
\begin{equation}\label{KK}
\KK _\lambda(g)=\int_s^\infty[\int_s^z\frac{G_\lambda^2(z)}{G^2_\lambda(u)}du]\frac{\gamma(\gamma+1)}{z^2}g(z)dz \ .
\end{equation}
For $s_0>0$ large enough, there exists a unique  $R_\lambda (s) \in L^\infty([s_0,\infty[)$ solution to 
\begin{equation}\label{KKbis}
(\mbox{Id}-\KK_\lambda)R_\lambda=\KK_\lambda(1).
\end{equation}
Moreover, $R_\lambda$ is holomorphic in $|\lambda|<\lambda_0$ and $R_\lambda(s)=O(s^{-\frac{3}{2}})$ uniformly in $|\lambda|<\lambda_0$.

\end{lemma}
\bp
We apply a fixed point theorem. First of all, there exists a constant $M$ such that for all $s\geq 1$, $|\lambda|\leq \lambda_0$ and $z\geq s$, we have 
$$
|\int_s^z\frac{G_\lambda^2(z)}{G^2_\lambda(u)}du|\leq \frac{M}{(1+|z|)^\frac{1}{2}}.
$$
Indeed,  let us denote  $U=\{(x+i\lambda)e^{i\frac{\pi}{6}},\quad x\geq 0, |\lambda|\leq \lambda_0\}$. For $z\in U, |z|\geq \frac{1}{2}$, 
$$Ai(z)= e^{-\frac{2}{3}z^\frac{3}{2}}\tau(z),\mbox{ with}\quad \frac{c_0}{(1+|z|)^\frac{1}{4}}\leq |\tau(z)|\leq \frac{c_1}{(1+|z|)^\frac{1}{4}}
$$
then 
$$\begin{array}{rcl}
\dps|\int_s^z\frac{G_\lambda^2(z)}{G_\lambda^2}du|&\leq& \dps C\int_s^z e^{-\frac{4}{3}\mbox{Re}([(z+i\lambda)^\frac{3}{2}-(u+i\lambda)^\frac{3}{2}]e^{i\frac{\pi}{4}})}du\\
&=&\dps Cz\int_{\frac{s}{z}}^1e^{-\frac{4}{3}z^{\frac{3}{2}}\mbox{Re}([(1+i\frac{\lambda}{z})^\frac{3}{2}-(t+i\frac{\lambda}{z})^\frac{3}{2}]e^{i\frac{\pi}{4}})}dt\\
&\leq & \dps Cz\int_0^\infty e^{-tz^\frac{3}{2}}dt\sim \frac{C}{z^\frac{1}{2}}\quad \mbox{if }z\geq s\geq 1.
\end{array}
$$
Thus for $|\lambda|\leq \lambda_0$
$$
|s^{n+\frac{3}{2}}\KK_\lambda(g)(s)|\leq M\gamma(\gamma+1)\int_s^\infty \frac{s^{n+\frac{3}{2}}}{z^{n+\frac{3}{2}+1}}|g(z)z^n|dz.
$$
Then 
$$
||s^{n+\frac{3}{2}}\KK_\lambda(g)||_{L^\infty([1,\infty[)}\leq \frac{M\gamma(\gamma+1)}{(n+\frac{3}{2})}||s^ng||_{L^\infty([1,\infty[)}.
$$
Finally, $\KK_\lambda$ is bounded in $L^\infty([s_0,\infty[)$ with
$$
||\KK_\lambda||_{L^\infty([s_0,\infty[)}\leq \frac{2M}{3}\gamma(\gamma+1) s_0^{-\frac{3}{2}}\leq \frac{1}{2}\quad \mbox{if }s_0\mbox{ big enough.}
$$ 
Then (\ref{KKbis}) has a unique solution in $L^\infty([s_0,\infty[)$, $R_\lambda(s)$, holomorphic in $|\lambda|\leq \lambda_0$ 
and, since 
$\KK_\lambda(1)=O(s^{-\frac{3}{2}})$, we get the following asymptotics $R_\lambda(s)=O(s^{-\frac{3}{2}})$.

Moreover, since $\KK_\lambda^{(n+1)}(1)=O(s^{-\frac{3(n+1)}{2}})$, the sum $\sum_0^\infty \KK_\lambda^{(n+1)}(1)$ converges and we can write the asymptotic expansion

$$
R_\lambda=\sum_0^\infty \KK_\lambda^{(n+1)}(1).
$$

\ep
Let us resume the proof of Proposition \ref{propH}.
\noindent Since $Ai(e^{i\frac{\pi}{6}}(s+i\lambda))(1+R_\lambda(s))$ is solution on $[s_0,\infty[$, it may be extended on $]0,\infty[$ by $\tilde H_\lambda(s)$ in an holomorphic way for $|\lambda|<\lambda_0$ that can thus be written
$$
\tilde H_\lambda =a(\lambda) s^{-\gamma}F_{+,\lambda}(s)+b(\lambda)s^{\gamma+1}F_{-,\lambda}(s)
$$
where $a(\lambda)$ and $b(\lambda)$ are holomorphic for  $|\lambda|<\lambda_0$.
\noindent It remains now to prove that $a(0)\neq0$. For that purpose, assume that $a(0)=0$. 
Since $\gamma>\frac{1}{2}$, we get that the solution of  (\ref{eqdiffsing}) $\tilde H_0(s)\in L^2(\RR^+)$, and $\tilde H_0(s)=O(s^{-\infty})$, then by integration by parts, since $\tilde H_0(0)=0$ and since because of the asymptotic behavior of the chosen Airy function, $|\tilde H_0(s)|\leq_{s\sim 0} Cs^{\gamma+1}$, we write
$$
\int_0^\infty |\tilde H_0'(s)|^2+\frac{\gamma(\gamma+1)}{s^2}|\tilde H_0(s)|^2+is|\tilde H_0(s)|^2ds=0
$$
that leads to $\tilde H_0=0$ which leads to a contradiction.
To end the proof, we just define $H_\lambda(s)= \frac{1}{a(\lambda)}\tilde H_\lambda(s)$.
\ep

\begin{lemma}\label{lem1}
For all $\lambda\in \CC$, $\vert\lambda\vert\leq \lambda_0$ and $s\in]0,\infty[$, 
we have $H_\lambda(s)\neq 0$.
\end{lemma}

\bp
If $H_\lambda(s_0)=0$, since $H^{(k)}_\lambda=_{s\sim\infty}O(s^{-\infty})$ for any derivative of order $k\in \NN$, as above, by an integration by parts, we get
$$\int_{s_0}^\infty |H'_\lambda(s)|^2+[is+\frac{\gamma(\gamma+1)}{s^2}]|H_\lambda(s)|^2ds=\lambda\int_{s_0}^\infty |H_\lambda(s)|^2ds.
$$
It leads to 
$$
\int_{s_0}^\infty s|H_\lambda(s)|^2ds=
\mbox{Im}\lambda\int_{s_0}^\infty|H_\lambda(s)|^2ds\quad \Rightarrow \mbox{Im}\lambda>0 \ ,
$$
$$
\int_{s_0}^\infty \frac{\gamma(\gamma+1)}{s^2}|H_\lambda(s)|^2ds\leq\mbox{Re}\lambda\int_{s_0}^\infty|H_\lambda(s)|^2ds\quad \Rightarrow \mbox{Re}\lambda>0.
$$
More precisely, by summing those two equations, we get
$$\mbox{Im}\lambda+\mbox{Re}\lambda\geq c_0=min_{s\geq 0}(s+\frac{\gamma(\gamma+1)}{s^2})$$  which contradicts the fact that $|\lambda|\leq \lambda_0$. Indeed, $c_0$ does not depend on $\lambda_0$,
and we can choose $\lambda_0$  as small as we want which leads to a contradiction.
\ep

In the remaining part of this section, we prove the fact that $d(0)\neq 0$. For that purpose, we prove that 
\begin{lemma}Let $F_{+,\lambda}$ be defined in (\ref{defdesF}), we have
$$
s^{-\gamma}F_{+,0}(s)\notin L^2([1,\infty[).
$$
\end{lemma}
\bp
{\bf {Step1 : Changing the unknown.}}

Following (\ref{lesgn}), we write
$$
F_{+,0}=\sum_0^\infty d_n s^{3n},\quad \mbox{with }d_0=1,\quad d_{n+1}=\frac{i}{9(n+1)(n+1-\alpha)}d_n,\mbox{ where }\alpha= \frac{2\gamma+1}{3}\notin \NN
$$
that can also be written
$$
F_{+,0}=D_\alpha(s^3),\quad \mbox{where } D_\alpha= \sum_0^\infty d_nx^n.
$$
By introducing the sequence $ h_n$ writing 
$$
d_n=(\frac{i}{9})^n\frac{h_n}{n!}
$$
we get the following recurrence formula
\begin{equation}\label{defhn}
h_0=1\quad \mbox{ and }h_{n+1}=\frac{h_n}{n+1-\alpha}.
\end{equation}
Note that it implies that
$$
|d_n|\leq 9^{-n} (\frac{1}{n!})^2
$$
which implies that 
$$
|D_\alpha(x)|\leq Ce^{C\sqrt x}.
$$
Define now 
$$
F_\alpha(z)= \sum_0^\infty h_n z^n,
$$
we have for all $t\in\CC$, $\mbox{Re}t>0$,
\begin{equation}\label{reldalphafalpha}
\int_0^\infty e^{-\frac{x}{t}}D_\alpha(x)dx=tF_\alpha(\frac{it}{9}).
\end{equation}
Let us now study $F_\alpha$.

{\bf{Step2 : Study of $F_\alpha$}}

\begin{lemma}Let $F_\alpha=\sum_0^\infty h_n z^n$ with $h_n$ satisfying (\ref{defhn}). One has
\begin{equation}
F_\alpha(z)=\Gamma(1-\alpha) z^\alpha e^z+\alpha \int_0^\infty \frac{e^{-zv}}{(1+v)^{\alpha+1}}dv\quad \forall z\in \CC \mbox{ such that }\mbox{Re}z>0.
\end{equation}
\end{lemma}
\bp
Since the sequence $h_n$ is defined by (\ref{defhn}), the function $F_\alpha$ satisfies the differential equation
$$
F'_\alpha -\frac{\alpha}{z}(F_\alpha-F_\alpha(0))=F_\alpha, \quad F_\alpha(0)=1.
$$
By integrating this equation, we get for $z>0$
$$
F_\alpha (z) =C(\alpha) z^\alpha e^z+z^\alpha e^z\int_z^\infty \frac{\alpha e^{-s}}{s^{\alpha+1}}ds.
$$
But by integration by part, we write
$$
\int_z^\infty \frac{\alpha e^{-s}}{s^{\alpha+1}}ds=\frac{1}{z^\alpha}e^{-z}+\frac{1}{1-\alpha}\int_z^\infty \frac{(\alpha-1) e^{-s}}{s^{\alpha}}ds.
$$
By iterating this process, we obtain
$$\begin{array}{rcl}
F_\alpha(z)&=&\dps1+\frac{z}{1-\alpha}+\cdots+\frac{z^n}{(1-\alpha)\cdots(n-\alpha)}+z^\alpha e^z[C(\alpha)+\frac{1}{(1-\alpha)\cdots(n-\alpha)}\int_z^\infty e^{-s}s^{n-\alpha}ds]\\
\\
&=&\dps1+\frac{z}{1-\alpha}+\cdots+\frac{z^n}{(1-\alpha)\cdots(n-\alpha)}\\
&&+z^\alpha e^z[C(\alpha)-\Gamma(1-\alpha)+\frac{\Gamma(1-\alpha)}{\Gamma(n-\alpha+1)}\int_0^z e^{-s}s^{n-\alpha}ds]
\end{array}
$$
By letting $n\rightarrow \infty$, we get that
$$
C(\alpha)= \Gamma(1-\alpha)
$$
where $\Gamma(z)=\dps \int_0^\infty x^{z-1}e^{-x}dx$. Which concludes the proof of the lemma.
\ep

\bigskip

{\bf {Step3 : Proof that $D_\alpha$ is not a tempered distribution}}

Going back to (\ref{reldalphafalpha}), we obtain
$$
\int_0^\infty e^{-\frac{x}{t}}D_\alpha(x)dx=tF_\alpha(i\frac{t}{9})=t[\Gamma(1-\alpha)(i\frac{t}{9})^\alpha e^{i\frac{t}{9}}+\alpha\int_0^\infty \frac{e^{-i\frac{t}{9v}}}{(1+v)^{\alpha+1}}dv]
$$
then the Laplace transform of $1_{x>0}D_\alpha(x)$ is given by
\begin{equation}\label{ftrans}
\int_0^\infty e^{-\lambda x}D_\alpha(x)dx=[\frac{\Gamma(1-\alpha)}{\lambda}(\frac{i}{9\lambda})^\alpha e^{i\frac{1}{9\lambda}}+9\alpha\int_0^\infty \frac{e^{-iw}}{(1+9\lambda w)^{\alpha+1}}dw].
\end{equation}
Since $1_{x>0} D_\alpha(x)\leq Ce^{C\sqrt x}$, the Fourier transform 
$$
\int_0^\infty e^{-ix\xi} 1_{x>0} D_\alpha(x)dx$$ exists and is holomorphic in $\mbox{Im} \xi<0$ and from (\ref{ftrans}), we get
\begin{equation}\label{fouriertrans}
\int_0^\infty e^{-ix\xi}D_\alpha(x)dx=\frac{\Gamma(1-\alpha)}{i\xi}(\frac{1}{9\xi})^\alpha e^{\frac{1}{9\xi}}+9\alpha \int_0^\infty \frac{e^{-iw}}{(1+9iw\xi)^{\alpha+1}}dw.
\end{equation}
In (\ref{fouriertrans}), the term $\frac{\Gamma(1-\alpha)}{i\xi}(\frac{1}{9\xi})^\alpha e^{\frac{1}{9\xi}}$ is not the Fourier transform of a tempered distribution, unlike the second term. Indeed
$$
9\alpha \int_0^\infty \frac{e^{-iw}}{(1+9iw\xi)^{\alpha+1}}dw=\mathcal F(1_{x>0}\frac{\kappa(\alpha)}{2i\pi}\int_0^\infty e^{-i\frac{9x}{t}-t}t^{\alpha-1}dt).
$$
with $\kappa(\alpha)= \int_{\gamma_0}e^z\frac{dz}{z^\alpha}$, where $\gamma_0$ is the contour in $\CC$
connecting $-\infty$ to $-\infty$ with on loop conterclockwise around $z=0$.  Thus $D_\alpha$ is not tempered and $F_{+,0}(s)=D_{\frac{(1+2\gamma)}{3}}(s^3)$ is not tempered either.
\ep

\begin{lemma}\label{lemd0}The $L^2$ solution of (\ref{eqdiffsing}) is given by
$$
H_\lambda(s)=s^{-\gamma}F_{+,\lambda}(s)+d(\lambda)F_{-,\lambda}(s)s^{\gamma+1}
$$
where $d(\lambda)$ is an holomorphic function for $|\lambda|\leq \lambda_0$ and $d(0)\neq0$.
Moreover 
\begin{equation}\label{gld0}
d(0)= -e^{i\frac{\pi}{2}\frac{2\gamma+1}{3}} 9^{-\frac{2\gamma+1}{3}} \frac{\Gamma(1-\frac{2\gamma+1}{3})}{\Gamma(1+\frac{2\gamma+1}{3})}.
\end{equation}
\end{lemma}

\bp
If $d(0)=0$, then $s^{-\gamma}F_{+,0}=H_0\in L^2(1,\infty])$ which is a contradiction.
We need to compute the asymptotic of the non tempered part of $F_{+,\lambda}$ and $F_{-,\lambda}$ in order to compute the value of $d(0)$, the unique value such that this non tempered part vanishes.

Let us introduce $W_\alpha$, the non tempered part of $D_\alpha$, {\it{i.e.}}
$$
W_\alpha(x)= \frac{\Gamma(1-\alpha)}{2\pi} \int_{\Im \xi<0} e^{ix\xi+\frac{1}{9\xi}}\left(\frac{1}{9\xi}\right)^\alpha\frac{d\xi}{i\xi}.
$$
To compute the asymptotic, we use the stationary phase method with the phase $\phi(\xi)= x\xi-i\frac{1}{9\xi}$. The critical point corresponding to the point where $\phi'$ equals zero is given by $3\xi_c= e^{-i\frac{\pi}{4}}x^{-\frac{1}{2}}$ and we get
$$
W_\alpha(x)\sim_{x\sim\infty} c_0\Gamma(1-\alpha)3^{-\alpha} x^{\frac{\alpha}{2}-1}e^{i\frac{\pi}{4}\alpha}e^{\tau x^\frac{1}{2}},
$$
where $\tau= e^{i\frac{\pi}{4}}+\frac{1}{3} e^{-i\frac{\pi}{4}}$.
So 
$$
F_\alpha(s)\sim_{s\sim\infty}c_0\Gamma(1-\alpha)3^{-\alpha} s^{\frac{3\alpha}{2}-3}e^{i\frac{\pi}{4}\alpha}e^{\tau s^\frac{3}{2}}
$$
and finally, introducing $\mu=\frac{2\gamma+1}{3}$,
$$
H_0(s)\sim_{s\sim\infty}c_0 e^{\tau s^\frac{1}{2}}s^{-3}[\Gamma(1-\mu) 3^{-\mu}s^{3\frac{\mu}{2}}s^{-\gamma}e^{i\frac{\pi}{4}\mu}+d(0)\Gamma(1+\mu) 3^{\mu}s^{3\frac{-\mu}{2}}s^{1+\gamma}e^{-i\frac{\pi}{4}\mu}]
$$
which, since $s^{3\frac{\mu}{2}}s^{-\gamma}=s^{\frac{-3\mu}{2}}s^{1+\gamma}$ implies that
$$
d(0)= -e^{i\frac{\pi}{2}\frac{2\gamma+1}{3}} 9^{-\frac{2\gamma+1}{3}} \frac{\Gamma(1-\frac{2\gamma+1}{3})}{\Gamma(1+\frac{2\gamma+1}{3})}.
$$
\ep

Going back to the starting variable, introducing $\eta=\eps k$, $0<\eta\leq \eta_0$, and recall that 
$$
\LL^0_\eta=-\partial^2_v +\frac{\gamma(\gamma+1)}{v^2}+i\eta v
$$
we get the following proposition

\begin{proposition}
For any $\mu\in\CC$, $|\mu|\leq \eta^{\frac{2}{3}}\lambda_0$, the function
\begin{equation}\label{deftheta}\begin{array}{rcl}
\Theta_{\lambda,\eta}(v)&=&v^{-\gamma}F_{+,\lambda}(\eta^{\frac{1}{3}}v)+d(\lambda)F_{-,\lambda}(\eta^{\frac{1}{3}}v)v^{\gamma+1}\eta^{\frac{2\gamma+1}{3}}\\
\\
&=&\eta^{\frac{\gamma}{3}}H_\lambda(\eta^{\frac{1}{3}}v)
\end{array}\end{equation}
spans the space of solution in $L^2([1,\infty[)$ of the equation 
\begin{equation}\label{eqL0}
(\LL^0_\eta -\mu)g=0, \quad \mbox{with }\mu= \lambda \eta^{\frac{2}{3}}
\end{equation}
\end{proposition}
\bp
Define $g\in L^2(v_0,\infty)$ solution of (\ref{eqL0}), defining $v=\eta^{-\frac{1}{3}}s$ and $\mu=\eta^\frac{2}{3}\lambda$, the function $\tilde g(s)=g(\eta^{-\frac{1}{3}}s)$ satisfies $\tilde g \in L^2(s_0, \infty)$ and 
$$
(-\partial^2_s+\frac{\gamma(\gamma+1)}{s^2}+is-\lambda)\tilde g=0
$$
which ends the proof.
\ep

\subsection{Back to the real equation}

We consider now the complete operator,
$$
\LL_\eps=\LL^0_\eps+N(v)
$$
where $N(v)=W(v)-\frac{\gamma(\gamma+1)}{v^2}\in O(\frac{1}{v^4})$.

The goal of this section is to prove the following proposition.

\begin{proposition}\label{propexistlargevel}
There exists $\lambda_0$, $\eta_0$, such that the equation
$$
\left\{\begin{array}{rcl}
(\LL_\eps-\lambda\eta^\frac{2}{3})J_{\lambda,\eta}(v)&=&0,\quad v\in[0,\infty[\\
J_{\lambda,\eta}(0)&=&1
\end{array}\right.
$$has a  continuous solution in $(\lambda,\eta,v)\in \{|\lambda|\leq \lambda_0\}\times\{0\leq\eta\leq\eta_0\}\times[0,\infty[$, holomorphic in $\lambda\in \{|\lambda|<\lambda_0\}$ and satisfying $\int_0^\infty|J_{\lambda,\eta}(v)|^2dv<\infty$. Moreover this solution is unique.
\end{proposition}

As in the previous section, we will look for solutions in $L^2([v_0,\infty[)$, close to $\Theta_{\lambda,\eta}$ when $v\rightarrow\infty$ by writing
$$
G_{\lambda,\eta}=\Theta_{\lambda,\eta}(1+R_{\lambda,\eta}),\quad \mbox{where }R_{\lambda,\eta}(v)\rightarrow_{v\rightarrow\infty} 0.
$$

This change of unknown leads to the following equation for $R_{\lambda, \eta}$
\begin{equation}\label{eqrlambda}
\left\{\begin{array}{rcl}
(\mbox{Id}-\KK_{\lambda,\eta})R_{\lambda,\eta}&=&\KK_{\lambda,\eta}(1)\\
\KK_{\lambda,\eta}(g)(v)&=&\dps\int_v^\infty(\int_v^w \frac{\Theta^2_{\lambda,\eta}(w)}{\Theta^2_{\lambda,\eta}(u)}du)N(w)g(w)dw.
\end{array}\right.
\end{equation}
Note that by Lemma \ref{lem1}, we are allowed to divide by $\Theta^2_{\lambda,\eta}(u)$.

Before proving the proposition, we start with a series of lemma in order to proceed to a fixed point argument.

\begin{lemma}There exists $C_0$ such that for all $0<v<w$, we have 
\begin{equation}\label{boundtheta}
|\int_v^w\frac{\Theta^2_{\lambda,\eta}(w)}{\Theta^2_{\lambda,\eta}(u)}du|\leq C_0w\quad \forall |\lambda|\leq \lambda_0,\quad \forall 0<\eta\leq \eta_0.
\end{equation}

\end{lemma}
\bp
Back to the definition (\ref{deftheta}) of $\Theta$ and by writing $v=\eta^{-\frac{1}{3}}a$, $w=\eta^{-\frac{1}{3}}b$ and $u=\eta^{-\frac{1}{3}}t$, (\ref{boundtheta}) is true if and only if
$$
|\int_a^b\frac{H^2_\lambda(b)}{H_\lambda^2(t)}dt\leq C_0b
$$
then it is sufficient to prove that
\begin{equation}\label{boundhbis}
|\int_0^b\frac{H^2_\lambda(b)}{H_\lambda^2(t)}dt\leq C_0b.
\end{equation}

It is true if $b$ is small since $t\rightarrow |H_\lambda(t)|$ is decreasing near $0$. It is true for $b\in [b_0,B_0]$ compact set included in $]0,\infty[$. Finally, for $b>B_0$, we use the asymptotic coming from the Airy function 
\begin{equation}\label{boundh}
|H_\lambda(s)|\sim_{s\sim\infty} \frac{C}{s^\frac{1}{4}}e^{-\frac{\sqrt2}{3}s^\frac{3}{2}}.
\end{equation}
\ep

\begin{remark}
Note that for $b$ small, (\ref{boundhbis}) is  sharp, but for $b$ large, one can get the better estimate $C_0b^{-1/2}$.
\end{remark}

\begin{lemma}\label{lemg1}
There exists a function $G_{\lambda,\eta}(v)$ solution to $(\LL_\eps-\lambda\eta^\frac{2}{3})G_{\lambda,\eta}(v)=0$ for $v\in[0,\infty[$ and $G_{\lambda,\eta}$ is continuous in $\eta\in[0,\eta_0]$, holomorphic in $\lambda\in \CC$, $|\lambda|<\lambda_0$, continuous  in $(\lambda,\eta,v)\in \{|\lambda|\leq \lambda_0\}\times\{0\leq\eta\leq\eta_0\}\times[0,\infty[$ and there exists
$v_0>0$ such that
$$
G_{\lambda,\eta}=\Theta_{\lambda,\eta}(1+R_{\lambda,\eta}),\quad \mbox{with }|R_{\lambda,\eta}(v)|\leq\frac{C}{v^2}, 
\quad \text{for all} \ v\geq v_0
$$
where $C$ does not depend on $(\lambda,\eta)\in \{|\lambda|\leq \lambda_0\}\times\{0\leq\eta\leq\eta_0\}$. \\
Moreover, $G_{\lambda,0}$ does not depend on $\lambda$.
\end{lemma}

\bp
Since $N(w)=_{w\rightarrow\infty}O(\frac{1}{w^4})$, we get as in the proof of Lemma \ref{deuxdeux}
$$
||v^{n+2}\KK_{\lambda,\eta}(g)||_{L^\infty([1,\infty[)}\leq \frac{C_0}{n+2}||v^n g||_{L^\infty([1,\infty[)}.
$$
Then, there exists $v_0>>1$ that does not depend on $|\lambda|\leq \lambda_0$ and  $0<\eta\leq \eta_0$ such that (\ref{eqrlambda}) has a unique solution $R_{\lambda,\eta}(v)\in L^\infty([v_0,\infty[)$ and we have
$$
R_{\lambda,\eta}=O(\frac{1}{v^2}) \quad \mbox{when }v\rightarrow\infty.
$$

Moreover, thanks to (\ref{eqrlambda}), $R_{\lambda,\eta}(v) $ is an holomorphic function in $\{\lambda\in \CC,|\lambda|\leq \lambda_0\}$ for all $0<\eta\leq \eta_0$ as well as $G_{\lambda,\eta}$ for $v\in[0,\infty[$ since 
$$
G_{\lambda,\eta}=\Theta_{\lambda,\eta}(1+R_{\lambda,\eta}) \quad \text{for all} \ v\geq v_0
$$
 and $G_{\lambda,\eta}$ satisfies the differential equation 
 $$
 (\LL_\eps-\lambda\eta^\frac{2}{3})G_{\lambda,\eta}=0,\quad \forall v\in\RR.
 $$
 
 Note also that $ G_{\lambda,\eta}$ may be extended to $\eta=0$   and $G_{\lambda,0}$ does not depend on $\lambda$ since $\Theta_{\lambda,0}(v)=v^{-\gamma}$ and $\KK_{\lambda,0}(g) =\int_v^\infty \frac{w}{2\gamma+1}(1-(\frac{v}{w})^{2\gamma+1})N(w) g(w) dw$ do not depend on $\lambda$. Thus we get $G_{\lambda,0}=G_{0,0}(v)=v^{-\gamma}(1+O(v^{-2}))$. Continuity follows from the fact that thanks to (\ref{deftheta}) and (\ref{boundtheta}), the function $0<v<w$, $\lambda\in \CC$, $|\lambda|\leq\lambda_0$, $0\leq\eta\leq\eta_0$
 $$
 \int_v^w\frac{\Theta^2_{\lambda,\eta}(w)}{\Theta^2_{\lambda,\eta}(u)}du
 $$
 is holomorphic in $\lambda$, continuous in $\eta\in[0,\eta_0]$ and bounded by $C_0w$ which implies that we can apply the Lebesgue Theorem.
\ep

\bp{\bf Proof of Proposition\ref{propexistlargevel}}
First of all, since $(-\partial^2_v+W)G_{0,0}=0$, we have $\partial^2_vG_{0,0}=O(v^{-(\gamma+2)})$ then $\partial_v G_{0,0}=O(v^{-(\gamma+1)})$ (since $G_{0,0}=O(v^{-\gamma})$). Assume that $G_{0,0}(0)=0$, then by integration by parts of the collision operator written as in (\ref{defQ}), we get
$$
\int_0^\infty F[(\frac{G_{0,0}}{F^\frac{1}{2}})']^2=0
$$
then $G_{0,0}=CF^\frac{1}{2}$ and since $F(0)\neq 0$, we get $C=0$, then $ G_{0,0}=0$ that contradicts the fact that $G_{0,0}\sim_{v\sim\infty}v^{-\gamma}$.

Then, for $\lambda_0,\eta_0$ small, and $|\lambda|\leq\lambda_0$, $0\leq \eta\leq \eta_0$, we have $G_{\lambda,\eta}(s)\neq0$ since $G_{0,0}(0)\neq0$ and $G_{\lambda,\eta}(0)$ is continuous with respect to $\lambda$, $\eta$. Then $J_{\lambda,\eta}=\frac{G_{\lambda,\eta}}{G_{\lambda,\eta}(0)}$ is well defined. Uniqueness comes from the results above for $\eta>0$. When $\eta=0$, we also have uniqueness since the only solution of 
$$
[-\partial^2_v+W]f=0, \quad f(0)=0, \quad f\in L^2 
$$
is $f=0$.

\ep

\begin{remark}
The function $M(v)=\frac{1}{(1+v^2)^{\gamma/2}}$ is the unique solution in $L^2([0,\infty[)$ of the equation 
$(-\partial^2_v+W)f=0$ wich satisfies $f\simeq v^{-\gamma}$ for $v\rightarrow \infty$.
Since in the proof of Lemma \ref{lemg1} we have shown  $G_{\lambda,0}(v)= G_{0,0}(v)=v^{-\gamma}(1+O(v^{-2}))$, we get
\begin{equation}\label{gl1}
G_{\lambda,0}(v)=G_{0,0}(v)=M(v)\ .
\end{equation}

\end{remark}

\begin{proposition}\label{Glambda} Properties of $G_{\lambda,\eta}$.

\begin{itemize}
\item There exists a constant $C_0$ such that $\forall v\geq0$, $|\lambda|<\lambda_0$,  $\eta \in [0,\eta_0]$
\begin{equation}\label{GinfM} |G_{\lambda,\eta}(v)|\leq C_0 M(v).\end{equation}

\item We have the following limit
\begin{equation}\label{lim}
\lim_{\eta\rightarrow 0^+}\int_0^\infty \eta^{\frac{1}{3}}vG_{\lambda,\eta}(v) M(v) dv=0.
\end{equation}

\item For all $\lambda$, $\forall v$, 
\begin{equation}\label{GtendM}
\lim_{\eta\rightarrow 0} G_{\lambda,\eta}(v)=M(v).
\end{equation}

\end{itemize}

\end{proposition}
\bp
Concerning the first point, for $v\geq v_0$, we use the fact that the function $s^{\gamma}H_\lambda(s)$ is bounded on $[0,\infty[$, uniformly in $\vert \lambda\vert\leq\lambda_0$, and we write, with $s=\eta^{1/3}v$, 
$$
|G_{\lambda,\eta}(v)|= |\Theta_{\lambda,\eta}(v)(1+R_{\lambda,\eta})(v)|\leq C |\Theta_{\lambda,\eta}(v)|
=Cv^{-\gamma}| s^{\gamma}H_\lambda(s)|\leq C'v^{-\gamma}\sim C' M(v) \ .
$$
For $v\in [0,v_0]$, it follows from the continuity of $G_{\lambda,\eta}$, $G_{0,0}=M$ and $\min_{v\in [0,v_0]}M(v)>0$.\\
To prove the limit of the second point, we cut the expression in the following way
$$\begin{array}{rcl}
|\int_0^\infty\eta^{\frac{1}{3}}vGM|&\leq &\dps |\int_0^{s_0\eta^{-\frac{1}{3}}}\eta^{\frac{1}{3}} vGM dv|+\eta^{\frac{2\gamma-1}{3}}\int _{s_0}^\infty s^{1-\gamma} |H_\lambda(s)|ds\\
\\
&\leq&\dps  C_0 s_0 \int_0^\infty  M^2+\eta^{\frac{2\gamma-1}{3}}\int _{s_0}^\infty s^{1-\gamma} |H_\lambda(s)|ds
\end{array}
$$
and we conclude by using $\gamma> 1/2$ and
$$
\lim_{\eta\rightarrow0}\eta^{\frac{2\gamma-1}{3}}\int _{s_0}^\infty s^{1-\gamma} |H_\lambda(s)|ds=0\quad  \forall s_0>0
$$
then, after passing to the limit in $\eta$, we pass to the limit when $s_0\rightarrow 0$.\\
The third point follows from \eqref{gl1} and the continuity with respect to   $\eta$ of $G_{\lambda, \eta}$.
\ep

\subsection{Computation of the eigenvalue}

In this subsection, we proceed to a reconnection of the two parts of the eigenvector, the positive velocity part and the negative velocity part. In order to be able to do the reconnection, we need to compute the derivative of the eigenvector at $v=0$. 

Let  $G_{\lambda, \mu}$ defined above and introduce  the notations
$$
a(\lambda,\eta)\mbox{ satisfying }a(\lambda, \eta){G_{\lambda,\eta}}_{|_{v=0}}=1,\quad \mbox{and}\quad b(\lambda,\eta)\mbox{ satisfying }a(\lambda, \eta){G'_{\lambda,\eta}}_{|_{v=0}}=b(\lambda,\eta)
$$
Observe that the functions $a(\lambda,\eta), b(\lambda,\eta)$ are holomorphic in 
$\lambda\in \CC, \vert\lambda\vert <\lambda_0$, and since $G_{\lambda,0}=M$, one has 
$a(\lambda,0)=1, b(\lambda, 0)=0$.\\
 Due to symetries in particular due to the parity of $M$, the connection condition reads $b(\lambda,\eta)+\overline b(\overline \lambda,\eta)=0$.
We thus need to compute $\Re b(0,\eta)$ and the coefficient in front of $\lambda$.
We gather all the needed results in the following proposition

\begin{proposition}\label{proplam0}
\begin{itemize}
\item The expression of $b(\lambda,\eta)$ is given by
\begin{equation}
b(\lambda,\eta) = a(\lambda,\eta) \eta^{\frac{2}{3}}\int_0^\infty (\lambda-i\eta^{\frac{1}{3}}v)G_{\lambda,\eta}(v) M(v)dv
\end{equation}
\item The coefficient in front of $\lambda$ is given by
\begin{equation}
\lim_{\eta\rightarrow 0^+} b(\lambda,\eta) \eta^{-\frac{2}{3}}=\lambda \int_0^\infty M^2(v)dv.
\end{equation}
\item concerning the real part of $b(0,\eta)$, we get
\begin{equation}
\lim_{\eta\rightarrow 0^+} \eta^{-\frac{2\gamma+1}{3}}\Re b(0,\eta)=\int_0^\infty s^{1-\gamma}\Im (H_0(s))ds=
(2\gamma+1)\Re (d(0))
\end{equation}
where 
$$
d(0)= -e^{i\frac{\pi}{2}\frac{2\gamma+1}{3}} 9^{-\frac{2\gamma+1}{3}} \frac{\Gamma(1-\frac{2\gamma+1}{3})}{\Gamma(1+\frac{2\gamma+1}{3})}.
$$
\end{itemize}

\end{proposition}
\bp
The first point is obtained by integrating the equation satisfied by $G_{\lambda,\eta}$ by part.\\
To get the second point, we use  $\lim_{\eta\rightarrow 0^+}G_{\lambda,\eta}=M$ which implies
$\lim_{\eta\rightarrow 0^+}a(\lambda,\eta)=1$, and we conclude by using \ref{lim}. 

The computation of $\Re b(0,\eta)$ will be split into three steps.
Recall that 
$$
b(0,\eta)= -i\eta a(0,\eta) \int_0^\infty wG_{0,\eta}(w) M(w)dw
$$
In order to get the result, we prove the three following lemmas. 
\begin{lemma}\label{petitesvitesses}

 The small velocities don't participate to the limit of the coefficient $b(0,\eta)$, 
\begin{equation}
\lim_{\eta\rightarrow 0^+} \eta^{-2\frac{(\gamma-1)}{3}}\int_0^{v_0} w\Im [a(0,\eta)G_{0,\eta}]Mdw=0
\end{equation}

\end{lemma}

\begin{lemma}\label{grandesvitessesImG} We have
\begin{equation}\label{eqsgammaH0}
\lim_{\eta\rightarrow 0^+} \eta^{-2\frac{(\gamma-1)}{3}}\int_{v_0}^\infty w\Im [a(0,\eta)G_{0,\eta}] M(w) dw=\int_0^\infty s^{1-\gamma} \Im H_0(s) ds
\end{equation}
\end{lemma}

In order to prove those results, we need the following lemma
\begin{lemma} \label{ima}For all $\gamma>1$, we have
 \begin{equation}\label{imG1}
 |\Re (aG)-M|\leq C\eta 
 \end{equation}
 
  \begin{equation}\label{imaG1}
 |\Im (aG)|\leq C\eta.  \end{equation}

Moreover, for large velocities,
 \begin{equation}\label{imG}
 |\Re (aG)-M|\leq C\eta <v>^{3-\gamma}, \quad \forall v \in [v_0, s_0\eta^{-\frac{1}{3}}]
 \end{equation}
 
  \begin{equation}\label{imaG}
 |\Im (aG)|\leq C\eta <v>^{3-\gamma}, \quad \forall v \in [v_0, s_0\eta^{-\frac{1}{3}}].
 \end{equation}
 
\end{lemma}

\bp[proof of lemma \ref{ima}]

Set $f_\eta = \Re (aG)$, and $\eta l_\eta= \Im (aG)$. They satisfy the following equations,
with $Q=-\partial^2+W$

\begin{equation}\label{eqfeta}
Q[f_\eta]-\eta^2vl_\eta=0,\quad f_\eta(0)=1
\end{equation}

\begin{equation}\label{eqleta}
Q[l_\eta]+vf_\eta=0,\quad l_\eta(0)=0
\end{equation}
By multiplying the equation by $M$ and integrating by parts,
we compute their derivatives
$$
f'_\eta(0)=\eta^2\int_0^\infty vl_\eta M\quad \mbox{and}\quad l'_\eta(0)=-\int_0^\infty vf_\eta M \ .
$$
Lemma\ref{ima} can be reformulated as follows 

\begin{equation}\label{boundfeta}
f_\eta(v)=M(v)+\tilde f_\eta\mbox{ with }|\tilde f_\eta |\leq C\eta <v>^{3-\gamma} \ ,
\end{equation}
and \begin{equation}\label{boundleta}
| l_\eta |\leq C <v>^{3-\gamma}.
\end{equation}

The solution of $Q(f)=g, \quad f(0)=a \mbox{ and }f'(0)=b$
is given by
$$
f=-\int_0^v g(w)M(w) dwZ(v)+\int_0^v g(w)Z(w) dwM(v)+aM(v)+ bZ(v)
$$
where $$M(v)= {(\frac{F}{C_\beta})}^\frac{1}{2}(v)=\dps{\frac{1}{(1+|v|^2)^{\frac{\gamma}{2}}}}\sim_{v\sim  \infty} v^{-\gamma} \quad \mbox{and}\quad\dps Z(v)= M(v) \int_0^v \frac{1}{M^2(w)}dw\sim_{v\sim  \infty} v^{\gamma+1}.$$

Since the function $f_\eta$ satisfies
$$
Q[f_\eta] =\eta^2vl_\eta,\quad f_\eta(0)=1, \quad f'_\eta(0)=\eta^2\int_0^\infty vl_\eta M.
$$
we get 
$$
f_\eta(v)=M(v)+ \eta^2(\int_0^\infty vl_\eta M)Z(v)+ (\int_0^v \eta^2 vl_\eta Z)M(v) -(\int_0^v \eta^2vl_\eta M)Z(v)
$$
which can be rewritten 
$$
f_\eta= M(v)+ \tilde f _\eta (v)
$$
where 
$$
\tilde f_\eta (v) =(\int_0^v \eta^2 vl_\eta Z)M(v) +(\int_v^\infty \eta^2vl_\eta M)Z(v).
$$
Since $|aG|\leq CM$, we get both $|f_\eta|\leq CM$ and $|\eta l_\eta|\leq CM$. Since $\gamma>1$, $vM^2$ is integrable at infinity and we write $\dps \int_v^\infty vM^2\leq C <v>^{2-2\gamma}$ and we finally get (\ref{boundfeta}).

Concerning $l_\eta$, it satisfies the equation
$$
Q[l_\eta]= -vf_\eta, \quad l_\eta(0)=0, \quad l'_\eta(0)=-\int_0^\infty vf_\eta M
$$
which leads to the following formula
$$
\begin{array}{rcl}
l_\eta(v)&=&\dps -(\int_0^\infty vf_\eta M) Z(v)-(\int_0^v vf_\eta Z) M(v)+(\int_0^v vf_\eta M) Z(v)\\
&=&\dps -(\int_v^\infty vf_\eta M) Z(v)-(\int_0^v vf_\eta Z) M(v).
\end{array}
$$
As before, since $\gamma>1$ and $f_\eta\leq CM$, we get (\ref{boundleta}).

\ep

\bp[Proof of Lemma \ref{petitesvitesses}]
Case 1 : $\gamma\in ]1,\frac{5}{2}]$.

First of all, since $2(\gamma-1)/3 <1$, and $\vert \Im[a(0,\eta)G_{0,\eta}]\vert =\vert \eta l_\eta\vert \leq C\eta$, we get
$$
\eta^{-2\frac{(\gamma-1)}{3}} w\Im [a(0,\eta)G_{0,\eta}(w)]M(w)\rightarrow_{\eta\rightarrow 0}0 \quad \text{for all} \ w.
$$
But since $|a(0,\eta)G_{0,\eta}|\leq CM$, when $\gamma>1$, one has 
$
w|a(0,\eta)G_{0,\eta}|M\leq CwM^2 \in L^1
$
and we conclude by the Lebesgue theorem.\\
Case 2: $\gamma\in[\frac{1}{2},1]$.
Since  $2(\gamma-1)/3 \leq 0$, we obtain directly the result by using the Lebesgue theorem and the third point of Proposition \ref{Glambda}  that gives 
$$
\int_0^{v_0} w a(0,\eta)G_{0,\eta}Mdw\rightarrow _{\eta\rightarrow 0} \int_0^{v_0} wM^2dw
$$
thus the imaginary part goes to zero.

\ep

\bp[Proof of Lemma \ref{grandesvitessesImG}]
 In order to prove (\ref{eqsgammaH0}), we proceed to a change of variable $w=\eta^{-\frac{1}{3}} s$,  which means that we need to compute
$$
\lim _{\eta\rightarrow 0^+} \int_{\eta^{\frac{1}{3}} v_0}^\infty \Im[a(\eta,0)\eta^{-\frac{\gamma}{3}}G_{0,\eta}(\eta^{-\frac{1}{3}}s)]s\eta^{-\frac{\gamma}{3}}M(\eta^{-\frac{1}{3}}s)ds
$$
where $\eta^{-\frac{\gamma}{3}}G_{0,\eta}(\eta^{-\frac{1}{3}}s)= H_0(s) [1+R_{0,\eta}(\eta^{-\frac{1}{3}}s)]$.
For that purpose, we use the Lebesgue Theorem, by writing that
$$
\forall s>0 \quad, \lim _{\eta\rightarrow 0^+}  \Im[a(\eta,0)\eta^{-\frac{\gamma}{3}}G_{0,\eta}(\eta^{-\frac{1}{3}}s)]s\eta^{-\frac{\gamma}{3}}M(\eta^{-\frac{1}{3}}s)= s^{1-\gamma} \Im (H_0(s)) .
$$
To obtain the domination, we use $\beta=2\gamma$ with $\beta\in ]1,5[\setminus \{2,3,4\}$.
Therefore, one has $\gamma \in ]1/2,1[\cup ]1,5/2[$. When $\gamma \in ]1,5/2[$, 
we use
$
M(w)= \frac{1}{(1+w^2)^\frac{\gamma}{2}}\leq |w|^{-\gamma},
$
which leads for $|s|>\eta^{\frac{1}{3}}v_0$, 
$$
|\Im[a(\eta,0)\eta^{-\frac{\gamma}{3}}G_{0,\eta}(\eta^{-\frac{1}{3}}s)]s\eta^{-\frac{\gamma}{3}}M(\eta^{-\frac{1}{3}}s)|\leq C|\Im [a(\eta,0)\eta^{-\frac{\gamma}{3}}G_{0,\eta}(\eta^{-\frac{1}{3}}s)]s^{1-\gamma} \ .
$$
Moreover, $\Im[aG]= \eta l_\eta$ and since $\gamma>1$, we have  for any $v\in [v_0, s_0\eta^{-1/3}]$, 
$|l_\eta|\leq C|v|^{3-\gamma}$.
So for $s\leq s_0$, we get
$$
\vert \Im[a(\eta,0)\eta^{-\frac{\gamma}{3}}G_{0,\eta}(\eta^{-\frac{1}{3}}s)]\vert \leq Cs^{3-\gamma}
$$
and since $\dps \gamma<\frac{5}{2}$
$$
\vert \Im[a(\eta,0)\eta^{-\frac{\gamma}{3}}G_{0,\eta}(\eta^{-\frac{1}{3}}s)]s\eta^{-\frac{\gamma}{3}}M(\eta^{-\frac{1}{3}}s)\vert \leq Cs^{4-2\gamma}\in L^2(]0,1])\ .
$$
For $s\geq s_0$, we use the fact that $|a(0,\eta)]\leq C$ and $|R_{0,\eta}|\leq C$ and we write
$$
\vert \Im[a(\eta,0)\eta^{-\frac{\gamma}{3}}G_{0,\eta}(\eta^{-\frac{1}{3}}s)]s\eta^{-\frac{\gamma}{3}}M(\eta^{-\frac{1}{3}}s)\vert
\leq C|H_0(s)|s^{1-\gamma}\in L^1[1,\infty[)
$$
since $\dps H_0(s)\sim_\infty s^{-\frac{1}{4}}e^{-\frac{\sqrt2}{3}s^\frac{3}{2}}$.
When $\gamma\in ]1/2,1[$, we just use $H_0(s)\sim_0 s^{-\gamma}$, and we write 
$$
\vert \Im[a(\eta,0)\eta^{-\frac{\gamma}{3}}G_{0,\eta}(\eta^{-\frac{1}{3}}s)]s\eta^{-\frac{\gamma}{3}}M(\eta^{-\frac{1}{3}}s)\vert
\leq C|H_0(s)|s^{1-\gamma}\in L^1]0,\infty[) \ .
$$
Then since the function is dominated by an integrable function, we can pass to the limit and we conclude that 
(\ref{eqsgammaH0}) holds true.
\ep

\begin{lemma}[Computation of the coefficient]
The coefficient of the leading power in $\eta$ of the real part of $b(0,\eta)$ given in Lemma \ref{grandesvitessesImG}
 is equal to 
\begin{equation}
\int_0^\infty s^{-\gamma} s \Im H_0 ds= (1+2\gamma) \Re d(0) .
\end{equation}
\end{lemma}

\bp Recall that $H_0$ satisfies 
$$
P(H_0)= -isH_0, \quad P(f)=(-\partial^2_s +\frac{\gamma(\gamma+1)}{s^2})f
$$
that implies $s\Im H_0 =\Re(P H_0)=P(\Re H_0)$.
In another hand, 
$$H_0(s)=s^{-\gamma} F_{+,0}(s)+d(0) s^{\gamma+1}F_{-,0}(s)=
s^{-\gamma} \big(1+\frac{is^3}{6(1-\gamma)}+O(s^6) + d(0) s^{2\gamma+1}(1+O(s^3))\big)\ .$$
Since $\gamma\in ]1/2,1[\cup]1,5/2[$, this  implies that $s^{-\gamma} s\Im H_0$ is integrable when $s\sim 0$.
Moreover, since $P(s^{-\gamma})=0$, we can proceed to a double integration by part by writing
$$
\int_0^\infty s^{-\gamma} s\Im H_0 ds=\lim_{s_0\rightarrow 0} \int_{s_0}^\infty s^{-\gamma}  P(\Re H_0) ds  
=\lim_{s_0\rightarrow 0} [\partial_s \Re H_0 s_0^{-\gamma} -\Re H_0 \partial_s (s^{-\gamma})]=(2\gamma+1)\Re d(0).
$$
Recall that $\Re d(0)\neq 0$ and the  computation of $d(0)$ has been done in Lemma\ref{lemd0}.
\ep
The proof of Proposition \ref{proplam0} is complete.
\ep

\subsection{Extension to the negative velocities and computation of the eigenvalue with lowest absolute value .}
Until this subsection, all the computations have been done for non negative velocities. We now need to extend this solution to negative velocities. For that purpose, we need to make a $C^1$ connection by connecting the value and the derivative at $v=0$.

\begin{proposition}\label{proplimmueps} 
Let $\eta_0>0$ and $\lambda_0>0$ small enough. For all $\eta\in [0,\eta_0]$, 
there exists in the complex disc $\{\mu\in \mathbb C, \vert\mu\vert\leq \eta^{2/3}\lambda_0\}$ 
a unique $\mu(\eta)$ such that the equation  (\ref{eqmeps}) (with $\eta=\varepsilon k$) admits a  
solution $M^\eta$  in $L^2(\RR)$. Moreover,  this solution is unique, and  one has

\begin{equation}\label{glres}
\begin{aligned}
&\mu(\eta) = \kappa \eta^{\frac{2\gamma+1}{3}}(1+O( \eta^{\frac{2\gamma+1}{3}}))\\
&\kappa= 2C_\beta^2 (2\gamma+1)9^{-\frac{2\gamma+1}{3}} \cos(\frac{\pi}{2}\frac{2\gamma+1}{3})
\frac{\Gamma(1-\frac{2\gamma+1}{3})}{\Gamma (1-\frac{2\gamma+1}{3})} >0\ .
\end{aligned}
\end{equation}
\end{proposition}

\bp Recall that the equation we consider is given by 

$$
-\partial^2_v +\frac{\gamma}{(1+|v|^2)^2}[|v|^2(\gamma+1)-1]+i\eta v)M^\eta=\mu M^\eta.
$$
If we change $v$ into $-v$, the equation remains the same except that we have to change $i$ into $-i$ (note that we assume here the parity of the equilibrium $M$) which means that 
$$
M^\eps_{\mu^\eps}(v)=\overline{M^\eps}_{\overline{\mu^\eps}}(-v).
$$
Thus, if we want to reconnect the derivative for $v=0$ in order to have a $C^1(\RR) $ function, we get the constraint
${M^\eps}'_{\mu^\eps}(0)=-{\overline{M^\eps}}'_{\overline{\mu^\eps}}(0)$ which is equivalent to
$$
b(\lambda, \eta)+\overline{b}(\overline\lambda,\eta)=0.
$$
By Proposition \ref{proplam0} and the normalization $C_\beta^2\int M^2 dv =1$, one has 
$\eta^{-2/3}b(\lambda, \eta)=b(0,\eta)+ \frac{\lambda}{2C_\beta^2}(1+o_\eta(1))+ O(\lambda^2)$, thus the connection equation reads
$ \Re (b(0,\eta))+ \frac{\lambda}{2C_\beta^2}(1+o_\eta(1))+ O(\lambda^2)=0$. This implies 
$\lambda=-2C_\beta^2\Re (b(0,\eta))+O((\Re (b(0,\eta)))^2)$. Then the result follows by the third point
in Proposition \ref{proplam0}, since by  Lemma \ref{lemd0},
formula \eqref{gld0}, one has $2C_\beta^2 (1+2\gamma)\Re (d(0))=-\kappa$.
\ep
\begin{remark}
For $\eta\in [-\eta_0,0]$, by complex conjugaison on the equation, we get 
$$\mu(\eta)=\overline{\mu(-\eta)}= \kappa \vert \eta\vert ^{\frac{2\gamma+1}{3}}(1+O(\vert  \eta\vert ^{\frac{2\gamma+1}{3}})).$$
\end{remark}

\section{Proof of Theorem \ref{main} :  Momentum method}
\label{secmoment}

\subsection{A priori estimates}

We start with a compactness Lemma.  
\begin{lemma}\cite{NaPu}\label{norm}
For initial datum $f_0 \in Y^p_\omega$ where $p\geq 2$ and a positive time $T$.
\begin{enumerate}
\item The solution $f^\varepsilon$ of \eqref{fp-theta} is bounded in $L^\infty\left([0, T]; \ Y^p_\omega )\right)$ uniformly with respect to $\varepsilon$ since it satisfies 
\begin{equation}\label{estimate1.1}||f^\varepsilon(T)||^p_{Y^p_\omega}+\frac{p\ (p-1)}{\theta(\varepsilon)}\int_0^T\int_{\mathbb{R}^{2d}}\frac{|\nabla_v(f^\varepsilon\ \omega)|^2}{\omega}\ (f^\varepsilon)^{p-2}\ \omega^{p-2}\ \mathrm{dv}\mathrm{dx}\mathrm{dt}\leq ||f_0||^p_{Y^p_\omega}.
\end{equation}
\item The density $\rho^\varepsilon(t,x)=\int_{\mathbb{R}^d}f^\varepsilon \ \mathrm{d}v$ is such that
\begin{equation} 
||\rho^\varepsilon(t)||^p_{p}\leq C_\beta^{-2(p-1)} ||f_0||^p_{Y^p_\omega}\quad \mathrm{for\ all}\quad t\in [0,T].
\end{equation}
\item Up to a subsequence, the density $\rho^\varepsilon$ converges weakly star in $L^\infty ([0, T]; L^p(\mathbb{R}^d))$ to $\rho$.
\item 
Up to a subsequence, the sequence $f^\varepsilon$ converges weakly star in $L^\infty ([0, T]; Y^p_\omega({\mathbb{R}^{2d}}))$ to $f=\rho(t,x)\frac{C_\beta^2}{\omega}.$
\end{enumerate}
\end{lemma}

\begin{corollary} Let $F=\frac{C_\beta^2}{\omega}=C_\beta^2 M^2$, $M=\frac{1}{(1+v^2)^\frac{\gamma}{2}}$.
Let $f^\eps$ solution to (\ref{fp-theta}) with $\theta(\varepsilon)=\varepsilon^{\frac{2\gamma+1}{3}}$.  Assume that  $||f_0\omega||_{\infty}\leq C$.
 Then  $g^\eps=f^\eps F^{-1/2}$ satifies the following estimate
\begin{equation}\label{estblan}
\int_0^T\int _{\mathbb R}\left(\int |g^\eps-\rho^\eps F^{1/2}|^2dv\right)^{\frac{2\gamma+1}{2\gamma-1}} ds dy\leq C\eps^\frac{2\gamma+1}{3}.
\end{equation}
\end{corollary}

\bp
Recall the Nash type inequality \cite{CGGR}\cite{rw} \cite{BaBaCaGu}:  for any $h$ such that $\int h Fdv=0$, we have
\begin{equation}\label{Nash}
\int h^2 Fdv\leq C\left(\int|\nabla_v h|^2 Fdv\right)^\frac{2\gamma-1}{2\gamma+1}(||h||_\infty^2)^\frac{2}{2\gamma+1} \ .
\end{equation}
Define $h=g^\eps F^{-1/2}-\rho^\eps =\frac{f^\eps}{F}-\rho^\eps$, define $\alpha=\frac{2\gamma+1}{3}$.
Observe that  from $||f||_{Y^p_\omega}=||\omega f||_{L^p(\frac{dxdv}{\omega})}$ 
and  Proposition \ref{norm}, formula \eqref{estimate1.1}, we have
$$
||h_0||_{L^\infty}=\lim_{p\rightarrow\infty } ||h_0||_{Y^p_\omega}\geq\lim_{p\rightarrow\infty } ||h||_{Y^p_\omega}
\geq||h||_{L^\infty} \ .
$$
 Thus by Lemma \ref{norm}, formula \eqref{estimate1.1}, we get 
$$\begin{array}{rcl}
\int_0^T\int _{\mathbb R} \left(\int |g^\eps-\rho^\eps F^{1/2}|^2dv\right)^{\frac{2\gamma+1}{2\gamma-1}} dsdy
&=& \int_0^T\int _{\mathbb R} \left(\int h^2 Fdv\right)^{\frac{2\gamma+1}{2\gamma-1}}dsdy \\ 
&\leq& 
C\int_0^T\int _{\mathbb R} \left(\int|\nabla_v h|^2 Fdv\right)(||h||_\infty^2)^\frac{2}{2\gamma-1} ds dy\\
&\leq &C\dps \int_0^T\int _{\mathbb R} \left(\int \frac{|\nabla_v(f^\eps \omega)|^2}{\omega}dv\right) ds dy
\leq C \eps^\alpha.
\end{array}$$

\ep

\subsection{Weak limit}

Recall $T=\varepsilon^{-\alpha}, \alpha=\frac{2\gamma+1}{3}$. By solving equation (\ref{rescaled}), we write
$$\tilde g^\eps(s,v, k)=e^{-sT\mathcal{L}_\eps}\tilde g(0,v,k)
$$
which gives going back to the rescaled space variable $y$
$$
 g^\eps(s,v, y)=\frac{1}{2\pi} \int e^{iy\cdot k}\tilde g^\eps(s,v,k)dk \ .
$$
Our purpose is to pass to the limit when $\eps\rightarrow 0$, or $T\rightarrow \infty$.\\
Recall $f^\eps(s,y,v)\geq 0$ and $\int f^\eps(s,y,v)dxdv=\int f_0(x,v)dxdv$ for all $s\geq 0$.\\
Let $\hat \rho^\eps (s,k)= \int e^{-iyk}\rho^\eps (s,y)dy$ be the Fourier transform in $y$ of 
$\rho^\eps = \int f^\eps dv= \int g^\eps F^{1/2} dv$ .

\begin{proposition}\label{lem:equilibrium}
 For  all $k\in \mathbb R$,  $\hat\rho^\varepsilon(.,k)$ converges  to $\hat \rho(.,k)$, 
unique solution  to the ode
\begin{equation}\label{fracdifeq}
\partial_s\hat\rho+\kappa|k|^\alpha \hat \rho=0, \quad \hat \rho_0=\int_\RR\hat  f_0 dv \ .
\end{equation}
\end{proposition}

\bp
Recall that $\mathcal{L}_\eps=Q+i\varepsilon k v$, $Q=-\partial_v^2+W$. Let $k\in \mathbb R$, $\eta=\varepsilon k$, and let
$M^\eta(v)$ be the unique solution in $L^2(\mathbb R)$ of $\mathcal{L}_\eps(M^\eta)=\mu(\eta)M^\eta$
given in Proposition \ref{proplimmueps}.  One has
$$
\begin{array}{rcl}
\frac{d}{ds}\int  \tilde g^\eps(s,v,k) M^\eta dv&=& \int \partial_s\tilde g^\eps M^\eta dv
= -\eps^{-\alpha}\int \mathcal{L}_\eps(\tilde g^\eps) M^\eta dv\\
&=&
-\eps^{-\alpha}\int \tilde g^\eps    \mathcal{L}_\eps(M^\eta) dv = -\eps^{-\alpha}\mu(\eta )\int \tilde g^\eps  M^\eta dv \ .  
\end{array}
$$
Therefore one has, with $F^\eps (s,y)=C_\beta \int   g^\eps(s,v,y) M^\eta dv$, 
\begin{equation}\label{gl9}
 \hat F^\eps (s,k) =e^{-s\eps^{-\alpha}\mu(\eps k)}  \hat F^\eps (0,k) \quad \forall s\geq 0.
\end{equation} 
By Proposition \ref{proplimmueps}, we have
$ \eps^{-\alpha}\mu(\eps k)\rightarrow \kappa \vert k\vert^\alpha$. 
Moreover, the following limit holds true:
\begin{equation}\label{gl11}
\forall k\in \mathbb R, \quad \hat F^\eps (0,k)=C_\beta\int \tilde g^\eps(0,v,k)  M^\eta dv  \rightarrow  \hat \rho_0(k) \ . 
\end{equation}
The verification of \eqref{gl11} is easy. One has $\tilde g^\eps(0,v,k)= \hat f_0(v,k) F^{-1/2}$
and $C_\beta F^{-1/2} M^\eta(v)=\frac{M^\eta}{M}(v) \rightarrow 1$ for all $v\in \mathbb R$ since 
our construction gives $M^\eta(v)=a(\lambda,\eta)G_{\lambda,\eta}(v)$ with  $a(\lambda,0)=1, \ G_{\lambda, 0}= M$.
Moreover, one has by \eqref{GinfM} the domination $\vert M^\eta(v)\vert \leq CM(v)$. Thus \eqref{gl11} holds true by Lebesgue Theorem.\\
\begin{remark}
Observe that it is only in the verification of \eqref{gl11} (initial data at time $s=0$) that we use the fact that 
$M^\eta$ is associated to the eigenvalue of smallest absolute value of the operator $ \mathcal{L}_\eps$,
since it is the only eigenfunction which satisfy $\lim_{\eta\rightarrow 0}M^\eta= M$. 
\end{remark}
It remains to verify 
\begin{equation}\label{gl10}
\forall k\in \mathbb R, \quad C_\beta\int \tilde g^\eps(s,v,k)  M^\eta dv  \rightarrow  \hat \rho(s,k)  
\quad \text{in} \ \mathcal D'(]0,\infty[).
\end{equation}
By \eqref{gl9} and \eqref{gl11}, for all $k\in \mathbb R$ and $s\geq 0$, one has
$\lim_{\eps\rightarrow 0}\hat F^\eps (s,k)=e^{-s\kappa \vert k\vert^\alpha}\hat\rho_0(k)$, thus \eqref{gl10}
will be consequence of the weaker
\begin{equation}\label{gl10bis}
\quad C_\beta\int  g^\eps(s,y,v)  M^\eta dv  \rightarrow   \rho(s,y)  
\quad \text{in} \ \mathcal D'(]0,\infty[\times \mathbb R) \ .
\end{equation}

\noindent Let us now verify \eqref{gl10bis}. For that purpose, we write
$$
C_\beta \int  g^\eps  M^\eta dv  -\rho= 
C_\beta  \int ( g^\eps- \rho^\eps F^{1/2})M^\eta dv +  \rho^\eps \int (C_\beta M^\eta-F^{1/2})F^{1/2} dv 
+\rho^\eps -\rho \ .
$$
By using (\ref{estblan}), (\ref{GinfM}) and (\ref{lim}), and the Lebesgue theorem we pass to the limit.
The proof of Proposition \ref{lem:equilibrium} is complete.
\ep

{\bf Proof of The main result: Theorem \ref{main}}. From the two last items in Lemma \ref{norm}, we have just to prove that the 
for any given $k$, the Fourier transform $\hat\rho(s,k)$ of the weak limit
$\rho(s,y)$, is solution  of the equation \eqref{diff}, which is precisely Proposition \ref{lem:equilibrium}.

\end{document}